\documentclass[11pt,a4paper,reqno]{amsart}
\usepackage[T1]{fontenc}
\usepackage{lmodern,amsmath,amssymb,mathtools,microtype}
\usepackage[margin=29mm]{geometry}
\usepackage[hidelinks]{hyperref}
\hypersetup{pdftitle={A simple proof of the Slicing Theorem via Gaussian transport},pdfauthor={Christian Kipp and Emanuel Milman}}
\newtheorem{theorem}{Theorem}[section]
\newtheorem{corollary}[theorem]{Corollary}
\newtheorem{lemma}[theorem]{Lemma}
\newtheorem{remark}[theorem]{Remark}
\newcommand{\R}{\mathbb R}
\newcommand{\E}{\mathbb E}
\newcommand{\Id}{\mathrm{Id}}
\newcommand{\dd}{\,\mathrm d}
\DeclareMathOperator{\tr}{tr}
\DeclareMathOperator{\Cov}{Cov}
\DeclareMathOperator{\Var}{Var}
\DeclareMathOperator{\diver}{div}
\newcommand{\HS}{\mathrm{HS}}
\newcommand{\cH}{\mathcal H}
\newcommand{\ip}[2]{\langle #1,#2\rangle}
\newcommand{\norm}[1]{\lVert #1\rVert}
\newcommand{\dto}{_{\#}}
\numberwithin{equation}{section}

\title[Slicing via Gaussian transport]{A simple proof of the Slicing Theorem\\via Gaussian transport}
\date{}

\author[Christian Kipp and Emanuel Milman]{Christian Kipp\textsuperscript{*,$\dagger$} and Emanuel Milman\textsuperscript{*,$\ddagger$}}

\thanks{$^\dagger$ Email: ckipp@campus.technion.ac.il. $^\ddagger$ Email: emilman@tx.technion.ac.il.}
\thanks{$^*$ Department of Mathematics, Technion -- Israel Institute of Technology, Haifa, Israel. \\
The research leading to these results is part of a project that has received funding from the European Research Council (ERC) under the European Union's Horizon 2020 research and innovation programme (grant agreement No 101001677).}

\begin{document}
\begin{abstract}
We give a simple short proof of Bourgain's Slicing Problem, resolved by Klartag and Lehec, by using optimal transport from the Gaussian measure. Our argument yields the (non-sharp) explicit bound $L_\mu < 1$ for the isotropic constant of any even log-concave measure $\mu$, and $L_K < \sqrt{2}/e$ for any convex body $K$. The proof is entirely self-contained, and does not use any tools from asymptotic geometric analysis such as (reverse) Blaschke--Santal\'o, (reverse) Brunn--Minkowski, the equivalence to small-ball estimates, a reduction from Thin-Shell, nor the existence of the moment map, Eldan's Stochastic Localization or Klartag's refined Lichnerowicz estimate. The idea is the same recurring theme as in the recent proofs of the sharp thin-shell estimate and the log-Brunn--Minkowski conjecture -- all the information is already contained in the second derivative of the change-of-variables formula of optimal transport. While the resolutions of the latter two problems used the moment-map method, the idea to use a source Gaussian measure is borrowed from OpenAI's recent solution of the strong Slicing Problem, which is based on much more convoluted arguments and computations.  \end{abstract}
\maketitle

\section{Introduction}

An absolutely continuous probability measure $\mu=f(x) \dd x$ on $\R^n$ is called log-concave if $\log f : \R^n \rightarrow \R \cup \{-\infty\}$ is concave. 
Let $\Sigma_\mu=\Cov(\mu)$ denote its covariance matrix, and define its isotropic constant, differential entropy and affine entropy as follows:
\[
  L_\mu :=\norm f_\infty^{1/n}(\det\Sigma_\mu)^{1/(2n)},\quad
  h(\mu) :=-\int_{\R^n} f\log f\dd x , \quad \cH(\mu):=h(\mu)-\tfrac12\log\det\Sigma_\mu . 
\]
Both $L_\mu$ and $\cH(\mu)$ are invariant under invertible affine maps. When $\Sigma_\mu = \Id$, we have  $\cH(\mu)=h(\mu)$. It is well-known (see e.g. \cite[Proposition I.2]{BobkovMadiman}) that 
\[
-\log \norm{f}_\infty \leq h(\mu) \leq n - \log \norm{f}_\infty ,
\]
or equivalently
\begin{equation} \label{eq:intro-HL}
-\log L_\mu \leq \frac{1}{n} \cH(\mu) \leq 1 - \log L_\mu . 
\end{equation}
When $\mu = \mu_K$ is uniform on a convex body $K$, one defines $L_K := L_{\mu_K}$ and trivially verifies that both inequalities on the left are equalities.

\medskip

Bourgain's Slicing Problem \cite{BourgainOrigin} asks whether $L_K$ is bounded by an absolute constant, uniformly over all dimensions and all convex bodies (see also \cite{MP89} for equivalent formulations). The strong Slicing Problem further predicts that $L_K \leq L_{\Delta_n}$, where $\Delta_n$ is a simplex in $\R^n$; in the class of centrally-symmetric convex bodies $K$, it predicts that $L_K \leq L_{Q_n}$, where $Q_n$ is the cube in $\R^n$. It was shown by Fradelizi and Mar\'in Sola \cite{FradeliziMarinSola} that, when formulated simultaneously in all dimensions, these two sharp conjectural upper bounds are equivalent to the corresponding sharp entropy lower bounds $\cH(\mu) \geq n$ and $\cH(\mu) \geq \frac{n}{2}\log 12$ for general and even
log-concave probability measures on $\R^n$ (attained by products of exponential
measures and by $\mu_{Q_n}$), respectively.
In dimension $n=1$, these sharp entropy bounds were established by
Melbourne, Nayar and Roberto \cite[Theorem~1.1]{MelbourneNayarRoberto} and Madiman, Nayar and Tkocz
\cite[Theorem~1, $p=2$]{MadimanNayarTkocz}, respectively. 

\medskip
 
After decades of contributions by numerous authors \cite{Bourgain,Klartag,Chen,KlartagLehec,Lich,Guan}, Klartag and Lehec finally gave a positive answer to the Slicing Problem in \cite{KL} using Eldan's Stochastic Localization \cite{Eldan}. A simplified proof was soon after obtained by Bizeul \cite{Bizeul}. A deterministic proof (without Stochastic Localization) was recently obtained by Brazitikos \cite{Brazitikos}, using only standard tools in asymptotic geometric analysis, as well as Klartag's refined Lichnerowicz estimate \cite{Lich}. It was shown by Eldan and Klartag \cite{EK} that a universal thin-shell estimate would imply a positive resolution of the Slicing Problem (see also the recent alternative argument by Brazitikos and Milman \cite{BM}); such a universal thin-shell bound was finally obtained by Klartag and Lehec \cite{KLthin} using Stochastic Localization, and recently in sharp form by Chen and Klartag \cite{CK} using the moment-map method. Just a few days ago, Klartag and Lehec found yet another explicit reduction from the sharp anisotropic thin-shell estimate using Stochastic Localization, which yields an explicit (non-sharp) bound $L_K \leq \sqrt{5 / \pi}$ \cite{KLnote}. On October 6, 2026, OpenAI announced the resolution of the strong Slicing Problem for general (possibly non-symmetric) convex bodies, obtaining the sharp bound
\[
 L_K \leq L_{\Delta_n} = \frac{(n!)^{1/n}}{(n+1)^{(n+1)/(2n)}\sqrt{n+2}} \nearrow \frac{1}{e} \;\;\; \text{(as $n \rightarrow \infty$)}.
\]
However, the OpenAI manuscript is very difficult to digest and comprehend, overwhelming the reader with a plethora of parameters and extensive numerical computations, whose purpose often remains unclear. 
The present work arose from a (failed) attempt to properly understand and simplify their argument. 

\medskip

In this work, we give a simple proof of the Slicing Theorem with an explicit constant. The proof is short and entirely self-contained, avoiding any tools from asymptotic geometric analysis such as (reverse) Blaschke--Santal\'o, (reverse) Brunn--Minkowski, the equivalence to small-ball estimates, a reduction from Thin-Shell, the existence of the moment map, Eldan's Stochastic Localization and Klartag's refined Lichnerowicz estimate. For this reason, we hope it may be of value and interest to the community.

\begin{theorem}[Even Case]\label{thm:main}
For every even log-concave probability measure $\mu$ on $\R^n$,
\begin{equation}\label{eq:entropy}
 \cH(\mu)\ge c_e n,\qquad
 c_e :=\frac12\log\frac{15\pi e}{8}-\frac1{\sqrt{10}} \simeq 1.07044 >1 ,
\end{equation}
and consequently $L_\mu\le e^{1-c_e}<1$. For the uniform measure on a centrally symmetric convex body $K$ one has $L_K\le e^{-c_e} < \frac{1}{e}$.
\end{theorem}

Recall that the conjectured sharp constant in (\ref{eq:entropy}) is $c_e = \frac{1}{2} \log 12 \simeq 1.2425$ (attained by $\mu_{Q_n}$). To the best of our  knowledge, the above estimates improve over the best known numeric bounds in the even case (see also Section \ref{sec:conclude} for further improvements). For the purposes of establishing the Slicing Theorem, there is no real loss in generality in restricting to even densities --- the reduction of the general case to the even one is well-known \cite[Remark~2]{isomorphic}, since $L_\mu \leq \sqrt{2} L_{\mu \ast \mu_-}$ where $\mu_-(A) = \mu(-A)$. By passing to the Keith Ball body $B$ associated to $\textbf{1}_K \ast \textbf{1}_{-K}$, a delicate computation performed in \cite[Theorem~1.1]{Martin} similarly shows that $L_K \leq \sqrt{2} L_B$.  This establishes:

\begin{corollary}[General Case]\label{cor:general}
Every log-concave probability measure $\mu$ on $\R^n$ satisfies $L_\mu \leq \sqrt{2} e^{1-c_e}< \sqrt{2}$. 
For every convex body $K$ we have $L_K \leq \sqrt{2} e^{-c_e} < \frac{\sqrt{2}}{e}$.
\end{corollary}

If one prefers a fully self-contained proof which does not rely on the (elementary) reduction to the even case, we remark that evenness is not essential for our argument. A variant for a general log-concave measure $\mu$ yields
\begin{equation} \label{eq:ca0}
 \cH(\mu)\ge c^0_a n,\qquad
 c^0_a :=\left(\frac12\log\frac{8\pi e}{5}-\frac1{\sqrt2}\right) \simeq 0.60026 ,
\end{equation}
which implies $L_\mu \leq e^{1 - c_a^0} \simeq 1.4914$; we give the short verification at the end of Section~\ref{sec:end}. However, this yields a worse estimate than the one from Corollary \ref{cor:general}, which can be shown to imply (see \cite[Lemma 5.1, formula (5.1)]{FradeliziMarinSola})
\[
\cH(\mu)\ge c_a n,\qquad c_a := c_e - \log \sqrt{2} \simeq 0.72387 . 
\]
By OpenAI's \cite[Theorem~1.2]{entropy}, the sharp constant above is $c_a = 1$ (attained by a product of exponential measures). Theorem \ref{thm:main} gives a simple proof that every even log-concave measure $\mu$ satisfies this estimate (as $c_e > 1$). 

\medskip

The idea underlying the proof of Theorem \ref{thm:main} is the same recurring theme as in the recent proofs of the sharp thin-shell estimate by Chen and Klartag \cite{CK} and the log-Brunn--Minkowski conjecture by OpenAI \cite{logBM} -- all the information is already contained in the second derivative of the change-of-variables formula of optimal transport. While the resolutions of the latter two problems used the moment-map method, we borrow the idea to use a source Gaussian measure  from OpenAI's  resolution of the strong Slicing Problem \cite{entropy}. 

\medskip

The rest of this work is organized as follows. In Section \ref{sec:2} we derive the Gaussian transport identity. In Section \ref{sec:3}, we test it against two test functions. In Section \ref{sec:4}, we establish Theorem \ref{thm:main}. In Section \ref{sec:5} we comment on how to improve the constants $c_e,c_a$ further. 

\bigskip

\textbf{AI Declaration.} This approach was discovered after several weeks of using ChatGPT to better understand OpenAI's proof of strong Slicing (and failing). Motivated by the approach which worked so well for sharp Thin-Shell and log-Brunn--Minkowski, the authors asked ChatGPT to use the second derivative of the change-of-variables formula for the Gaussian transport as its starting point. Specifically, the authors asked ChatGPT to multiply this identity by functions $f(\tau)$ of the transport differential $\tau$ (explicitly $f(x) = x^k$ and $f(x) = \log(x)$), take the trace, integrate by parts, and use the favorable properties of Gaussian integration such as the Poincar\'e and log-Sobolev inequalities; this immediately resulted in the present argument. The authors were hoping to find a simple proof of strong Slicing, but have settled (for now) for an explicit (non-sharp) constant for Slicing via a simple self-contained argument. ChatGPT was used in drafting an initial version of this paper, which was extensively polished by the authors, who take full responsibility for its contents.

\section{Gaussian transport} \label{sec:2}

Throughout Sections~2--4 we assume $\mu=e^{-V(x)}\dd x$ with
\begin{equation}\label{eq:regular}
 V\in C^\infty(\R^n),\qquad
 \alpha\Id\preceq\nabla_x^2 V\preceq\beta\Id,\qquad 0<\alpha\le\beta<\infty .
\end{equation}
The passage to arbitrary log-concave densities will be given in Appendix~\ref{app:approx}. Condition \eqref{eq:regular} is preserved under invertible linear maps, with possibly different constants. 

Let $\gamma$ be the standard Gaussian measure on $\R^n$ with density $\rho$, and let $T=\nabla_z\varphi$ be the Brenier optimal-transport map for the quadratic cost from $\gamma$ to $\mu$. We use $z$ for the Gaussian variable and $x=T(z)$ for the target variable. Set $\tau=dT = \nabla_z^2\varphi$. Unless indicated otherwise, $\E$ denotes integration against $\gamma$.

We recall the classical transport facts used below. It was shown by Brenier \cite{Brenier} that the optimal coupling for squared Euclidean cost is unique and is induced by the gradient of a convex potential $\varphi$. Under \eqref{eq:regular}, interior regularity makes $T$ a smooth diffeomorphism \cite{CafReg,CorderoErausquinFigalli}. Caffarelli's contraction theorem \cite{CaffarelliContraction} applied to $T$ and to its inverse gives
\begin{equation} \label{eq:Caffarelli}
       \beta^{-1/2}\Id\preceq\tau(z)\preceq\alpha^{-1/2}\Id .
\end{equation}
Indeed the transport from a potential with Hessian at most $b\Id$ to one with Hessian at least $a\Id$ is $\sqrt{b/a}$-Lipschitz \cite[Theorem~2.2]{KolesnikovContraction}. 
We require the existence of these bounds for technical reasons, but their exact values are irrelevant for our entropy estimates.

\begin{lemma}[Good position]\label{lem:position}
For a probability density satisfying \eqref{eq:regular} there is a symmetric positive definite linear map $A$ such that, for the Brenier potential $\varphi_A$ from $\gamma$ to $A\dto\mu$ and its Hessian $\tau_A=\nabla_z^2\varphi_A$,
\begin{equation}\label{eq:position}
                    \E\tau_A=\Id .
\end{equation}
\end{lemma}
\begin{proof}
Consider the set of correlation matrices
\[
 \mathcal C=\big\{\E[X \otimes Z] : X\sim\mu,\ Z\sim\gamma
                       \text{ under a common coupling}\big\}.
\]
It is convex and compact. Indeed the couplings form a weakly compact set. To pass correlations to a weak limit, truncate to $|X|,|Z|\le R$. On $|X|>R$ the error in truncating $\E[X \otimes Z]$ is bounded by $(\E[|X|^2\mathbf{1}_{\{|X|>R\}}] \E|Z|^2)^{1/2}$, and the analogous bound holds on $|Z|>R$.
Both bounds depend only on the fixed marginals and tend
to zero as $R\to\infty$, so the truncation error tends to zero uniformly
over all couplings.

Let $T_0$ be the Brenier map to $\mu$ in its original position. Gaussian integration by parts gives $\E[T_0(Z) \otimes Z]=\E[dT_0(Z)]\succ0$. Thus the maximum of $\det M$ over $M \in \mathcal C$ is positive. Since the Gaussian is rotation invariant, $\mathcal C$ is closed under right orthogonal multiplication. Polar decomposition therefore allows us to choose a symmetric positive definite maximizer $M_*$. For every $M\in\mathcal C$, differentiating the determinant at $t=0$ along $M_*+t(M-M_*) \in \mathcal C$, for $0\le t\le1$, gives
\[
             \tr(M_*^{-1}M)\le \tr(M_*^{-1} M_*) = n .
\]

Now select $A=M_*^{-1}$. A coupling realizing $M_*$ attains the maximum $n$ of $\E\ip{AX}{Z}=\tr(AM)$, and therefore minimizes the quadratic transport cost $\E|AX-Z|^2$ from $Z$ to $AX$. By uniqueness, $AX=T(Z)$ for the Brenier map $T$ pushing forward $\gamma$ to $A\dto\mu$. A final Gaussian integration by parts gives
\[
      \E[dT(Z)]=\E[T(Z) \otimes Z]=A M_*=\Id .\qedhere
\]
\end{proof}

The change-of-variables formula for the Gaussian transport  reads
\begin{equation}\label{eq:ma}
 V(T(z))=\frac{|z|^2}{2}+\frac n2\log(2\pi)+\log\det\tau(z).
\end{equation}
Differentiating once gives
\begin{equation}\label{eq:first}
            \tau ((\nabla_x V) \circ T) =z+\nabla_z\log\det\tau .
\end{equation}

Equip $\R^n$ with the Riemannian metric $\tau_z(u,v)=\ip{\tau(z)u}{v}$ (somewhat abusing notation) and reference measure $\gamma$. 
Write $L=-\Delta_{\tau,\gamma}$ for the negative of the corresponding weighted Laplacian on $(\R^n , \tau , \gamma)$, namely,
\begin{equation} \label{eq:L}
 Lu=-\rho^{-1}\diver_z(\rho\tau^{-1}\nabla_z u)
    =-\tr(\tau^{-1}\nabla_z^2 u)+\ip{(\nabla_x V) \circ T}{\nabla_z u}.
\end{equation}
To see the second equality, differentiate the inverse matrix, use the symmetry of $\nabla_z^3\varphi$ to get $\diver_z(\tau^{-1})=-\tau^{-1}\nabla_z\log\det\tau$, and apply \eqref{eq:first}. For smooth $u,v$ with $v$ compactly supported, integration by parts gives
\begin{equation} \label{eq:ibp}
\E[(Lu)v]=\E\ip{\tau^{-1}\nabla_z u}{\nabla_z v}=\E\big[\tau(\operatorname{grad}_\tau u,\operatorname{grad}_\tau v)\big]. 
\end{equation}
Consequently $L$ is symmetric and nonnegative on $C_c^\infty(\R^n)\subset L^2(\gamma)$. Moreover, (\ref{eq:L}) gives $\E L v = 0$ for all $v \in C_c^\infty(\R^n)$. 
We apply $L$ to matrices entry by entry in fixed Euclidean coordinates.

\begin{lemma}[Transport identity]\label{lem:identity}
On $(\R^n,\tau,\gamma)$ the following transport identity holds:
\begin{equation}\label{eq:identity}
       \Id-L\tau=P+Q, 
\end{equation}
where $P$ and $Q$ are defined as
\begin{equation} \label{eq:KQ}
      \qquad P :=\tau ((\nabla_x^2 V) \circ T)\tau, \qquad  Q_{ij} :=\tr\!\big(\tau^{-1}(\partial_{z_i}\tau)\tau^{-1}(\partial_{z_j}\tau)\big).
\end{equation}
Both $P$ and $Q$ are positive semi-definite.
\end{lemma}
\begin{proof}
Apply $\partial_{z_i}\partial_{z_j}$ to \eqref{eq:ma}. Since $\partial_{z_i}\partial_{z_j}T_a=\partial_{z_a}\tau_{ij}$, we have 
\[
 \partial_{z_i}\partial_{z_j} (V \circ T) = P_{ij}+\ip{(\nabla_x V) \circ T}{\nabla_z\tau_{ij}} . 
 \]
On the right side,
\[
 \partial_{z_i}\partial_{z_j}\log\det\tau
   =\tr(\tau^{-1}\partial_{z_i}\partial_{z_j}\tau)-Q_{ij}
   =\tr(\tau^{-1}\nabla_z^2\tau_{ij})-Q_{ij}.
\]
Collecting terms and recalling the definition (\ref{eq:L}) of $L$ gives \eqref{eq:identity}. Convexity gives $P \succeq0$, and $Q$ is the Gram matrix of $\tau^{-1/2}(\partial_{z_i}\tau)\tau^{-1/2}$ for the Hilbert--Schmidt inner product and hence positive semi-definite. 
\end{proof}
\begin{remark}
Analogous identities for the moment-map transport have appeared in \cite{CK, logBM}, and the one above for the Gaussian transport also appears in \cite{entropy}. Early versions of these computations appeared already in the work of Kolesnikov \cite{Kolesnikov}. 
\end{remark}

\section{Two tests of the transport identity}\label{sec:tests} \label{sec:3}

In this section we keep the computations formal, so as not to occlude their essence; the rigorous justification of all integrations and integrations by parts below by a standard truncation argument is deferred to Appendix~\ref{app:cut}.

At any fixed point choose orthonormal coordinates in which $\tau$ has eigenvalues $\lambda_1,\ldots,\lambda_n>0$ on the diagonal. Write $C_{ijk} =\partial_{z_i}\tau_{jk}$ in these coordinates. The tensor $C$ is symmetric in all three indices. Every formula below uses a fixed orthogonal change of coordinates at the point, so no eigenvector field is differentiated.  For a matrix field $B$, we denote $|\nabla_z B|^2=\sum_k\norm{\partial_{z_k}B}_{\HS}^2$.

\subsection*{The trace test and logarithm estimate}
For $a,b>0$ put $\ell[a,b] :=(\log a-\log b)/(a-b)$, with $\ell[a,a] :=1/a$. In the chosen coordinates, the divided-difference formula for matrix derivatives \cite[Theorem~V.3.3]{Bhatia1997} gives
\[
 (\partial_{z_i}\log\tau)_{jk}=\ell[\lambda_j,\lambda_k]C_{ijk}.
\]
To quickly verify this, note that by the spectral theorem $\log \tau = \int_0^\infty \left ( (1 + s)^{-1} \Id  - (\tau+s\Id)^{-1} \right ) \dd s$, and hence its derivative is $\int_0^\infty(\tau+s\Id)^{-1} \partial_{z_i} \tau (\tau+s\Id)^{-1}\dd s$.  Consequently, $C_{ijk} = \partial_{z_i} \tau_{jk}$ is multiplied by $\int_0^\infty (\lambda_j+s)^{-1}  (\lambda_k+s)^{-1}\dd s = \ell[\lambda_j,\lambda_k]$. 
Applying Cauchy--Schwarz, we may estimate
\[
 \ell[a,b]^2
 =\left(\int_0^1\frac{\dd t}{(1-t)a+tb}\right)^2
 \le\int_0^1\frac{\dd t}{((1-t)a+tb)^2}=\frac1{ab}.
\]
Consequently $|\nabla_z\log\tau|^2\le\sum_{ijk}C_{ijk}^2/(\lambda_j \lambda_k)=\tr Q$ pointwise. But since $\E L \tau = 0$, taking the trace of \eqref{eq:identity} and integrating gives
\begin{equation}\label{eq:trace}
                 \E\tr Q=n-\E\tr P\le n ,
\end{equation}
and we conclude that
\begin{equation}\label{eq:log}
                      \E |\nabla_z\log\tau|^2\le \E \tr Q \le n .
\end{equation}

\subsection*{The quadratic test}
Multiply \eqref{eq:identity} by $\tau^2$, take the trace, and integrate. Since $\tr(\tau^2 P)\ge0$,
\[
 \E\tr\tau^2\ \ge\ \E\tr(\tau^2 Q)+\E\tr(\tau^2 L\tau).
\]
In our coordinates the first integrand is $\sum_{ijk} \lambda_i^2 C_{ijk}^2/(\lambda_j \lambda_k)$. To compute the second integral use \eqref{eq:ibp}, observing that
$(\partial_{z_k}\tau^2)_{ij}=(\lambda_i+\lambda_j)C_{kij}$. It follows that
\[
 \E\tr\tau^2\ \ge\
 \E\sum_{ijk}\left(\frac{\lambda_i^2}{\lambda_j \lambda_k}+\frac{\lambda_i+\lambda_j}{\lambda_k}\right)C_{ijk}^2 .
\]
The symmetry of $C_{ijk}^2$ allows each coefficient to be replaced by its cyclic average. Writing $(a,b,c)=(\lambda_i,\lambda_j,\lambda_k)$, the average of $a^2/(bc)$, $b^2/(ca)$ and $c^2/(ab)$ is at least one, since their product is one. The cyclic sum of $(a+b)/c$ consists of the three reciprocal pairs $a/b+b/a$, $b/c+c/b$ and $c/a+a/c$, so its average is at least two. Altogether,
\begin{equation}\label{eq:quad}
                    3\E |\nabla_z\tau|^2\le\E\tr(\tau^2) .
\end{equation}

\section{Bounding the entropy and the isotropic constant}\label{sec:end} \label{sec:4}

Integrating the change-of-variables formula \eqref{eq:ma} gives
\[
 h(\mu)=h(\gamma)+\E\log\det\tau
       = \tfrac n2\log(2\pi e) +\E\log\det\tau.
\]
Recalling that 
\begin{equation} \label{eq:cH}
\cH(\mu) = h(\mu) -\tfrac12\log\det\Sigma_\mu = \tfrac n2\log(2\pi e) +\E\log\det\tau -\tfrac12\log\det\Sigma_\mu
\end{equation}
is affine-invariant, we have the freedom to choose a position of $\mu$ when trading off the required lower bound on $\E \log \det \tau$ with the upper bound on $\log\det\Sigma_\mu$. Instead of naively using a position with $\Sigma_{\mu} = \Id$ to eliminate the latter, the position from Lemma \ref{lem:position} allows us to efficiently control the former while keeping tabs on the latter. 

For $u\in H^1(\gamma)$ we use the well-known Gaussian Poincar\'e inequalities
\[
 \Var_\gamma u\le\E|\nabla_z u|^2,\qquad
 \Var_\gamma u\le\tfrac12\E|\nabla_z u|^2\quad\text{if $u$ is even}.
\]
Indeed in the orthonormal Hermite polynomial basis write $u=\sum_\alpha c_\alpha H_\alpha$. Then $\Var_\gamma u=\sum_{|\alpha|\ge1}|c_\alpha|^2$ and $\E|\nabla_z u|^2=\sum_\alpha |\alpha|\,|c_\alpha|^2$, by Gaussian integration by parts. For even $u$, all terms of odd total degree vanish. These identities first hold for polynomials and extend by density.

\subsection{The even case}
Assume in this subsection that $\mu$ is even. The map $z\mapsto-T(-z)$ is another convex-gradient transport from $\gamma$ to $\mu$. Uniqueness gives $T(-z)=-T(z)$. Hence $\tau$ and $\log\tau$ are even functions of $z$.

Our choice of position from Lemma \ref{lem:position} preserves evenness and ensures that $\E \tau = \Id$. Set 
\[
v :=\E\norm{\tau-\Id}_{\HS}^2 = \E \tr (\tau^2) -  n . 
\]
Applying the even Poincar\'e inequality entrywise to $\tau$ and using \eqref{eq:quad} gives
\[
       v\le\tfrac12\E|\nabla_z\tau|^2\le \tfrac16 \E \tr(\tau^2) = \tfrac16(n+v),
       \quad \text{ and hence } \quad  v\le\frac n5.
\]
Applied to $\log\tau$, the even Poincar\'e inequality and \eqref{eq:log} give
\[
       \E\norm{\log\tau-\E\log\tau}_{\HS}^2\le \frac12  \E |\nabla_z\log\tau|^2 \leq \frac n2 .
\]

To control $\E \log\det\tau$  from these two variances, observe that $r\log r\ge r-1$ gives
$\E\tr(\tau\log\tau)\ge\E\tr(\tau-\Id)=0$. Thus, using $\E\tau=\Id$ and applying Cauchy--Schwarz,
\begin{align}
\notag  -\E\log\det\tau & = -\E\tr\log\tau \le \E\tr(\tau\log\tau)-\E\tr\log\tau\\
\label{eq:logdettau}  &=\E\tr\!\big[(\tau-\Id)(\log\tau-\E\log\tau)\big]
 \le \sqrt{\frac n5}\sqrt{\frac n2}=\frac n{\sqrt{10}} .
\end{align}

It remains to bound $\log\det\Sigma_\mu$ from above. Let $Z\sim\gamma$ and put $S(z)=T(z)-z$. The field $S$ is odd, and Gaussian integration by parts gives 
\begin{equation} \label{eq:orthogonality}
 \E[S(Z) \otimes Z]=\E(dS) = \E(\tau-\Id)=0.
\end{equation}
 Thus each component of $S$ has Hermite expansion supported on degrees at least three, so
\[
 \E|S(Z)|^2\le\tfrac13\E\norm{\tau-\Id}_{\HS}^2
 =\frac v3\le\frac n{15}.
\]
By orthogonality (\ref{eq:orthogonality}), 
\[
\tr\Sigma_\mu= \E|T(Z)|^2 = \E|S(Z)+Z|^2 = \E|S(Z)|^2 + n\le \tfrac{16}{15} n.
\]
Consequently, by the AM-GM inequality applied to the eigenvalues of $\Sigma_\mu$,
\begin{equation} \label{eq:logdetSigma}
 \log\det\Sigma_\mu\le  n\log\tfrac{16}{15} .
 \end{equation}
 Plugging (\ref{eq:logdettau}) and (\ref{eq:logdetSigma}) into (\ref{eq:cH}) gives the asserted \eqref{eq:entropy}:
\[
\cH(\mu) \geq c_e n , \quad 
 c_e = \frac12\log(2\pi e)-\frac1{\sqrt{10}} -\frac12\log\frac{16}{15}
     = \frac12\log\frac{15\pi e}{8}-\frac1{\sqrt{10}} \simeq 1.07044  .
\]

\begin{remark} \label{rem:better}
A slightly simpler way to bound $\tr\Sigma_\mu$ is to apply the ordinary Gaussian Poincar\'e inequality to the components of $T$. Since $dT=\tau$, this gives
\[
 \tr\Sigma_\mu=\sum_{i=1}^n \Var_\gamma T_i\le\E\tr(\tau^2)
     =n+v\le\tfrac65 n ,
\]
which would ultimately yield the slightly smaller constant
\[
 c_e^0 = \frac12\log(2\pi e)-\frac1{\sqrt{10}} -\frac12\log\frac{6}{5}
     = \frac12\log\frac{5\pi e}{3}-\frac1{\sqrt{10}} \simeq 1.01155 .
\]
\end{remark}

This establishes \eqref{eq:entropy} under \eqref{eq:regular}, and Appendix~\ref{app:approx} extends this to general even log-concave densities. 
Recalling (\ref{eq:intro-HL}) and the subsequent comments, we have $\log L_\mu \leq 1 - \frac{1}{n} \cH(\mu) \leq 1 - c_e$ and $\log L_K = -\frac{1}{n} \cH(\mu_K) \leq - c_e$, concluding the proof of Theorem~\ref{thm:main}. 

\medskip

For completeness, let us recall the simple argument used to obtain the right-hand-side of (\ref{eq:intro-HL}), which is particularly self-contained in the even case. By evenness of the log-concave density $f$, $f$ attains its maximum at the origin, so we may write $f=\norm f_\infty e^{-W}$ with $W : \R^n \rightarrow [0,\infty]$ convex and $W(0)=0$. For $t\ge1$, convexity gives $tW(x/t)\le W(x)$. Integration yields $\int e^{-t W}dx \ge t^{-n}\int e^{-W} dx$. Subtracting the two sides from $\int e^{-W} dx$, dividing by $t-1$, and letting $t\downarrow1$ gives $\int W f dx\le n$. Therefore
\[
   \log \norm f_\infty +   h(\mu)= \int W f dx\le n,
 \]
 and hence
 \[
     \log L_\mu = \frac{1}{n} \left (\log \norm f_\infty + h(\mu) - \cH(\mu) \right ) \le 1-\frac{\cH(\mu)}n .
\]

\subsection{The general case}

Let us also briefly mention how to handle the general non-even case (even though Corollary \ref{cor:general} ultimately yields a better bound). 

Translate $\mu$ to have mean zero and use the linear position from Lemma \ref{lem:position} which ensures $\E \tau = \Id$. By ordinary (non-even) Gaussian Poincar\'e and \eqref{eq:quad},  $v\le(n+v)/3$, so $v\le n/2$. Similarly \eqref{eq:log} gives $\E\norm{\log\tau-\E\log\tau}_{\HS}^2\le n$. As a result, the same covariance estimate for $\tau$ and $\log\tau$ yields $-\E\log\det\tau\le n/\sqrt2$. 

The field $S(z)=T(z)-z$ is no longer odd, but we still have $\E S = 0$ and orthogonality (\ref{eq:orthogonality}). Consequently every component of $S$ has Hermite expansion supported on degrees at least two, and hence 
\[
\E|S|^2\le\tfrac12\E\norm{dS}_{\HS}^2=\tfrac12 \E \norm{\tau - \Id}_{\HS}^2 = \frac{v}{2} \leq \frac{n}{4} . 
\]
 Orthogonality \eqref{eq:orthogonality} then gives $\tr\Sigma_\mu= \E |S|^2 + n \le \frac{5}{4} n$. The rest of the argument remains unchanged, yielding \eqref{eq:ca0}.

\section{Concluding Remarks} \label{sec:conclude} \label{sec:5}
 
The constants $c_e$ and $c_a$ we obtained in Theorem \ref{thm:main} and following Corollary \ref{cor:general} are definitely not the best possible one can obtain using our simple approach, as all of the bounds we used in Section \ref{sec:4} are lossy, especially the loose $\E \tr(\tau \log \tau) \geq 0$. By testing additional test functions $f(\tau)$ against the transport identity \eqref{eq:identity}, a more detailed analysis can push these constants further: for example, we are able to obtain $c_e \simeq 1.22251$ and $c_a \simeq 0.97650$ by using a quadratic polynomial $f$, applying Gaussian Poincar\'e to a linear combination $h(\tau)$ of $\log(\tau)$ and  $\tau-\Id$, choosing a position for which $\E[ a \tau + b h(\tau) + c \Id] = 0$, and numerically optimizing over the various coefficients. With some additional computer-assisted numerical optimization, it is possible to further push this to $c_e \simeq 1.22825$ and $c_a \simeq 0.99327$. 
However, we do not know how to obtain the conjectured $c_e = \frac{1}{2} \log 12$ or recover OpenAI's $c_a = 1$, so we have chosen to present the simplest argument that yields $c_e > 1$.

\appendix
\section{Justification of the calculations}

In the appendix we justify some of the formal computations used above. The arguments are standard, albeit somewhat tedious. 

\subsection{Truncation of the integration by parts}\label{app:cut}
Assume \eqref{eq:regular}, and recall by (\ref{eq:Caffarelli}) that $m\Id\preceq\tau\preceq M\Id$ with $0< m < M <\infty$. Recalling the definition of $Q$ in (\ref{eq:KQ}), we have $\tr Q = |\tau^{-1/2} (\nabla_z \tau) \tau^{-1/2}|^2$, and hence
\[
       |\nabla_z\tau|^2\le M^2\tr Q,\qquad
       |\tau^{-1}\nabla_z\tr\tau|\le\frac{\sqrt n\,M}{m}\sqrt{\tr Q}.
\]
We first establish global integrability. Choose smooth $0\le\eta_R\le1$, equal to one on $B_R$, supported on $B_{2R}$ and satisfying $|\nabla_z\eta_R|\le 2/R$. Test the trace of \eqref{eq:identity} by $\eta_R^2$ and integrate using \eqref{eq:ibp}. Positivity of $P$ gives
\begin{align*}
 \E[\eta_R^2\tr Q]
 &\le n+2\E\!\left[\eta_R|\nabla_z\eta_R|
                    \,|\tau^{-1}\nabla_z\tr\tau|\right] \le n+\tfrac12\E[\eta_R^2\tr Q]
             +\frac{2nM^2}{m^2}\E|\nabla_z\eta_R|^2 .
\end{align*}
Fatou's lemma yields $\E\tr Q<\infty$, so $\nabla_z\tau\in L^2(\gamma)$. 

Testing the trace of \eqref{eq:identity} once more against $\eta_R^2$
and integrating by parts gives
\[
 \E[\eta_R^2\tr Q]
 = n\E[\eta_R^2]-\E[\eta_R^2\tr P]
 -2\E\!\left[\eta_R
 \left\langle\tau^{-1}\nabla_z\tr\tau,\nabla_z\eta_R\right\rangle\right].
\]
By Cauchy--Schwarz and the pointwise bound above, the last term is
bounded in absolute value by
$\frac{4\sqrt n\,M}{mR}(\E\tr Q)^{1/2}$, which tends to zero as $R \rightarrow \infty$. 
Since $\tr Q$ is integrable and $0\preceq P\preceq\beta M^2\Id$,
dominated convergence gives $\E\tr Q=n-\E\tr P \leq n$,
proving \eqref{eq:trace}. Integrating the pointwise estimate
$|\nabla_z\log\tau|^2\leq\tr Q$ from Section~\ref{sec:tests} now gives
\eqref{eq:log}. As $m\Id\preceq\tau\preceq M\Id$, the field
$\log\tau$ is bounded and hence belongs to $H^1(\gamma)$ by \eqref{eq:log}.

For the quadratic test, trace \eqref{eq:identity} against $\eta_R^2\tau^2$. Applying \eqref{eq:ibp} entrywise
gives
\begin{align*}
 \E[\eta_R^2\tr(\tau^2 L\tau)]
 &= \E\sum_{i,j}\eta_R^2
    \left\langle\tau^{-1}\nabla_z\tau_{ij},
    \nabla_z(\tau^2)_{ij}\right\rangle  + 2\E\sum_{i,j}\eta_R(\tau^2)_{ij}
    \left\langle\tau^{-1}\nabla_z\tau_{ij},\nabla_z\eta_R\right\rangle .
\end{align*}
Since $\|\tau^2\|_{\mathrm{HS}}\leq\sqrt n\,M^2$ and
$\tau^{-1}\preceq m^{-1}\Id$, Cauchy--Schwarz bounds the second term
in absolute value by
\[
 \frac{2\sqrt n\,M^2}{m}
 (\E|\nabla_z\tau|^2)^{1/2}(\E|\nabla_z\eta_R|^2)^{1/2}
 \longrightarrow 0 .
\]
For the first term, the product rule gives
$|\nabla_z(\tau^2)|\leq 2M|\nabla_z\tau|$, so its integrand is bounded
in absolute value by $(2M/m)|\nabla_z\tau|^2$, which is integrable.
Finally, $\tr(\tau^2 P)$ is bounded and
$0\leq\tr(\tau^2Q)\leq M^2\tr Q$. We may therefore let $R\to\infty$
by dominated convergence in the tested identity, obtaining the
quadratic test of Section~\ref{sec:tests} and hence \eqref{eq:quad}.

\subsection{Approximation of an arbitrary log-concave density}\label{app:approx}

Let $f=e^{-V}$ be a general log-concave probability density, with $V : \R^n \rightarrow \R \cup \{+\infty\}$. 
It is well-known and easy to see that integrability implies that $f(x) \leq C e^{-c |x|}$ for some $c , C > 0$. 

Approximating $f$ with a smooth log-concave $f_\delta$ satisfying \eqref{eq:regular} is standard, see for instance \cite[Proof of Lemma~2.1, p.~14, (5.1)]{Lich}; for completeness, we describe the construction. 
Let $\gamma_\delta$ denote the centered Gaussian density with covariance
$\delta\Id$ and set
\[
    f_\delta(x) = e^{-V_\delta(x)} 
    := \frac{1}{Z_\delta}\,e^{-\delta |x|^2/2}(f*\gamma_\delta)(x) ,
\]
where $Z_\delta > 0$ normalizes the integral to one. Then $V_\delta$ is smooth and
\[
    \delta\Id \leq \nabla^2 V_\delta(x)
    \leq (\delta+\delta^{-1})\Id
    \qquad \text{for all } x\in\R^n.
\]
Thus $\mu_\delta=f_\delta\,dx$ satisfies \eqref{eq:regular}.
If $f$ is even, then so is $f_\delta$.

As $\delta\downarrow0$, $f*\gamma_\delta\to f$ almost everywhere and in
$L^1$. Hence $Z_\delta\to1$ and $f_\delta\to f$ almost everywhere and
in $L^1$. Since $f(x)\leq C e^{-c|x|}$, we also have, for all sufficiently
small $\delta$,
\[
    f_\delta(x)
    \leq \frac{C}{Z_\delta} e^{-c|x|}
         \int_{\R^n} e^{c|y|}\gamma_\delta(y)\,dy
    \leq C' e^{-c|x|}.
\]
Dominated convergence gives $\Cov(\mu_\delta)\to\Cov(\mu)$ and
$h(\mu_\delta)\to h(\mu)$, using $|t\log t|\leq C''\sqrt{t}$ for
$0\leq t\leq C'$. Finally, almost-everywhere convergence gives
$\|f\|_\infty\leq\liminf_{\delta\downarrow0}\|f_\delta\|_\infty$,
while the defining formula gives
$\|f_\delta\|_\infty\leq Z_\delta^{-1}\|f\|_\infty$.
Thus $\|f_\delta\|_\infty\to\|f\|_\infty$, and therefore $\cH(\mu_\delta) \rightarrow \cH(\mu)$ and $L_{\mu_\delta}\rightarrow L_\mu$. 
Consequently, the stated bounds for these parameters apply to general and even log-concave densities.

\end{document}